\documentclass[11pt,a4paper,notitlepage,final,onecolumn,oneside]{article}
\usepackage{exscale}
\usepackage{latexsym}
\usepackage{calc}
\usepackage{amssymb}
\usepackage{url}
\usepackage[dvips]{graphicx}

\newcommand{\proof}{{\bfseries Proof:\quad}}

\begin{document}

\newtheorem{lemma}{Lemma}
\newtheorem{proposition}[lemma]{Proposition}
\newtheorem{theorem}[lemma]{Theorem}
\newtheorem{definition}[lemma]{Definition}
\newtheorem{convention}{Convention}
\newtheorem{hypothesis}[lemma]{Hypothesis}
\newtheorem{conjecture}[lemma]{Conjecture}
\newtheorem{remark}[lemma]{Remark}
\newtheorem{example}[lemma]{Example}
\newtheorem{property}[lemma]{Property}
\newtheorem{corollary}[lemma]{Corollary}
\newtheorem{algorithm}[lemma]{Algorithm}


\title{The graph isomorphism problem is polynomial}
\author{Aleksandr Golubchik\\
        \small Bramscher Str. 57, 49088 Osnabrueck, Germany,\\
        \small e-mail: agolubchik@gmx.de}

\date{June 1, 2006}
\maketitle


\begin{abstract}
It is known that a graph isomorphism testing algorithm is polynomially
equivalent to a detecting of a graph non-trivial automorphism algorithm.
The polynomiality of the latter algorithm, is obtained by consideration of
symmetry properties of regular $k$-partitions that, on one hand,
generalize automorphic $k$-partitions (=systems of $k$-orbits of
permutation groups), and, on other hand, schemes of relations (strongly
regular $2$-partitions or regular $3$-partitions), that are a subject of
the algebraic combinatorics.

It is shown that the stabilization of a graph by quadrangles detects the
triviality of the graph automorphism group. The result is obtained by
lineariation of the algebraic combinatorics.\bigskip

\end{abstract}


\section{Introduction}

It is known that the graph isomorphism problem is equivalent by complexity
to the problem of exposure of orbits of a graph automorphism group, and
these two problems are equivalent to the problem of detecting of graph
non-trivial automorphism. This equivalence was considered by R. Marthon
\cite{marthon} (s. also \cite{koebler}) and independently by author
\cite{gol}. In given text we show that the latter problem is polynomial.

Attempts to find the complexity of graph isomorphism problem were gone in
two directions: group theory and computational theory, and till today no
way brought a result. The literature to the first way one can find in
\cite{klin} and to the second in \cite{koebler}

We will go the first way and, more exactly, study symmetry properties of
combinatorial objects that follows from symmetry properties of $k$-orbits.
Earlier this direction led to a sequence of graph stabilization algorithms
(Weisfeiler-Lehman algorithm and its generalizations), to notions of
strongly regular graph and distance regular graph, and also to developing
of algebraic combinatorics (that has different origins \cite{bannai}).

What was of principal in studying of graph isomorphism it was a
simplification of a model. If to consider a graph as a partition $L_2$ of
a Cartesian square $V^2$ of a finite set $V$ (i.e. a color digraph with
colored vertices), then it is evident that the less is $|L_2|$ (the coarse
is $L_2$) the graph is simpler; so all specialists were concentrated on
the consideration of simple graphs, where $|L_2|=3$, classes are
symmetrical, i.e.  $\langle v_1,v_2\rangle L_2 \langle v_2,v_1\rangle$ for
any $v_1,v_2\in V$, and one class is a diagonal ($\{\langle
v_i,v_i\rangle\}$, non-colored vertices). But in this case classes of
actual examples are very large and not observable.  About 1986 author
discovered for him that there exists another way of the problem
simplification, it is the case where $|L_2|$ is as large as possible and
therefore classes of $L_2$ are small.  This way led to construction that
gives a very simple local representation of difficulties of the problem
and shows the way of combinatorial problem solution.

The main achievement in combinatorics for last 30 years is a developing of
algebraic combinatorics that studies associative schemes of relations:
partitions of $V^2$ possessing certain symmetry properties. By studying
of schemes of relations it were obtained some important examples of
strongly regular graphs satisfying to (so-called) $k$-condition for $k>3$
\cite{klin}. But nevertheless it was not developed a conceptional theory
of symmetries of $k$-partitions that could open a new view on graph
isomorphism. From here appears an other idea (1983): to consider the
symmetry properties of $k$-partitions for any $k\geq 2$.

And the third algebraic idea (1983) came from consideration of
stabilization algorithm for a $3$-partition $L_3$. Pure intuitively
it is clear that a stabilization of $L_3$ is an equation on $n^3$
($n=|V|$) variables. But if $L_3$ is obtained from an initial partition
$L_2$, then we have the same equation on $n^2$ variables, so it have to
exist some overdetermination (in a linear equation system the number of
equations is greater than number of variables). The attempts to find a
corresponding system of algebraic equations were without success many
years and it was decided to search for a pure combinatorial solution,
because the applying of this direction was always successful.

So this text was planned as realization of a combinatorial solution. But
on the way of working on the text and thank to the text \cite{klin} and
the book \cite{bannai} suddenly it was found an algebraic approach that
was simpler as combinatorial. In that connection in this paper is given
the algebraic solution of graph isomorphism problem and the combinatorial
one will be represent in separate text later.


\section{$k$-partitions}

Let $V$ be a $n$-element set, $V^k$ be Cartesian power of $V$ and
$V^{(k)}$ be the non-diagonal part of $V^k$, i.e. any $k$-tuple from
$V^{(k)}$ consists of $k$ different coordinates. Under
\emph{$k$-partition} $L_k$ below we understand a partition of $V^{(k)}$.

We shall consider a $k$-partition $L_k$ with a set of \emph{coloring
functions} $F=\{f:L_k\rightarrow Cr\}$, where $Cr$ is a set of ``colors''
which can have different identity: numbers, vectors, tensors and other.
Under an automorphism group $Aut(L_k)$ we understand the maximal
permutation group $G(V)$ that maintains each class of $L_k$ and so each
function $f\in F$.

Let $P,Q$ be partitions of a set $M$, then $P\sqcup Q$ and $P\sqcap
Q$ denote the union and intersection of $P$ and $Q$. If $P$ is a
subpartition of $Q$, then we write $P\sqsubset Q$.

The action of $Aut(L_k)$ on $V^{(k)}$ forms a partition
$Orb_k(Aut(L_k))\sqsubset L_k$ that consists of orbits of this
action or (as one say) of \emph{$k$-orbits} of $Aut(L_k)$. If
$L_k=Orb_k(Aut(L_k))$, then we say that $L_k$ is an \emph{automorphic
partition} (or system of $k$-orbits of a permutation group $G=Aut(L_k)$).

It is convenient to represent a $k$-set ($k$-relation, $k$-class)
$U_k\subset V^{(k)}$ as a matrix $M(U_k)$, whose line is a $k$-tuple of
$U_k$ and $i$-th column consists of values of $i$-th coordinate of
$k$-tuples. So the matrix of a $k$-set is defined accurate to line order.
If $L_k$ is an automorphic partition, then its class is an
\emph{automorphic $k$-set} ($k$-orbit). An automorphic $k$-partition and
its classes possess evident symmetry properties. Consider those properties.


\section{Regular and $pq$-stable $k$-partitions}

We say that a $k$-partition $L_k$ is \emph{$s$-symmetrical}, if for any
class $U_k$ of $L_k$ any $k$-relation $U_k'$ that differs from $U_k$ by
order of coordinates (order of columns in matrix $M(U_k)$) belongs to
$L_k$.

Let $\alpha_k=\langle v_1,\ldots,v_k\rangle$ be a $k$-tuple, $l<k$ and
$\alpha_l=\langle v_{i_1},\ldots,v_{i_l}\rangle$ be a $l$-tuple that is a
projection ($l$-projection) of $\alpha_k$ on a subspace $W=\{i_j:j\in
[1,l],i_j\in [1,k]\}$, then, using projecting operator $\hat{p}(W)$,
we write $\alpha_l=\hat{p}(W)\alpha_k$. The set of all $k$
$(k-1)$-projections of $\alpha_k$ we write as $\hat{p}\alpha_k$. Here the
projections are considered in natural order of coordinates determinate by
$\alpha_k$. The reverse to $\hat{p}$ operator $\hat{q}$ assembles $k$
$(k-1)$-projections $\hat{p}\alpha_k$ in $k$-tuple $\alpha_k=
\hat{q}\hat{p}\alpha_k$. From this definition follows the action of
projecting and assembling operators on $k$-relations and $k$-partitions.

For a $k$-relation $U_k$ is $U_k\subset\hat{q}\hat{p}U_k$. So we call
$U_k$ \emph{$l$-full} if $\hat{q}^{k-l}\hat{p}^{k-l}U_k=U_k$ (here
$\hat{p}^{k-l}U_k$ is a set of $k\choose l$ projections of $U_k$ on
$l$-subspaces). A $k$-partition $L_k$ is $l$-full, if all its classes are
$l$-full: $\hat{q}^{k-l}\hat{p}^{k-l}L_k=L_k$. If $L_k$ is $l$-full, then
it is evidently $(l+1)$-full.

We say that a $k$-partition $L_k$ is \emph{$p$-symmetrical}, if
$\hat{p}L_k$ is a $(k-1)$-partition, i.e. any two $(k-1)$-projections of
any two classes of $L_k$ are either equal or disjoint. It is clear that
$\hat{p}\hat{q}L_k\sqsubset L_k$. So we say that a $k$-partition $L_k$ is
\emph{$pq$-stable} if $L_k=\hat{p}\hat{q}L_k$.

\begin{proposition}\label{p^l-symm}
Let a $k$-partition $L_k$ be $p$-symmetrical, then a $(k-1)$-partition
$L_{k-1}=\hat{p}L_k$ is $p$-symmetrical too.

\end{proposition}
\proof
We consider the case $k=3$ (for $k>3$ the proof is similar). Let
$A=Cl(\langle 1,2,3\rangle)$ be a class of $L_3$ containing $3$-tuple
$\langle 1,2,3\rangle$, $B=Cl(\langle 1,4,5\rangle)\in L_3$ and
$C=Cl(\langle 1,2,4\rangle)\in L_3$, so that the intersection of
projections of $A$ and $B$ on a subspace $W=\langle 1\rangle$ is not
trivial: $\hat{p}(\langle 1\rangle)A\neq \hat{p}(\langle 1\rangle)B$. But
the inequality is not valid because of $\hat{p}(\langle 1\rangle)A=
\hat{p}(\langle 1\rangle)\hat{p}(\langle 1,2\rangle)A=\hat{p}(\langle
1\rangle)\hat{p}(\langle 1,2\rangle)C=\hat{p}(\langle
1\rangle)\hat{p}(\langle 1,4\rangle)C=\hat{p}(\langle
1\rangle)\hat{p}(\langle 1,4\rangle)B=\hat{p}(\langle 1\rangle)B$.
$\Box$\bigskip
\begin{proposition}\label{p->s-symm}
Let a $k$-partition $L_k$ be $p$-symmetrical and $k>2$, then it is
$s$-symmetrical.

\end{proposition}
\proof
It is sufficient to consider the case $k=3$ and assume that $L_3$ is a
$2$-full $3$-partition. Further we use the induction on $n$. For $n=3$ the
statement is easy verified. Let $n=4$, and $L_3$ be $p$-symmetrical
$2$-full $3$-partition. Let $L_1=\hat{p}^2L_3\neq \{V\}$, then this case is
reduced to cases $n<4$ and therefore the statement is correct. Let
$L_1=\{V\}=\{\{v_1,v_2,v_3,v_4\}\}$ and let $L_3'$ be a $3$-partition
obtained from $L_3$ by removal from it all $3$-tuples containing value of
a coordinate $v_4$, then $L_3'$ is also $p$-symmetrical $2$-full
$3$-partition for that statement is correct. From here it follows that the
statement is correct for $L_3$. The generalization on any $n$ is evident.

\begin{proposition}\label{pq=qp}
Let $L_k$ be a $pq$-stable $(k-1)$-full partition, then
$\hat{p}\hat{q}L_k=\hat{q}\hat{p}L_k=L_k$\,.

\end{proposition}

Let $U_k$ be a $k$-relation and $mU_l$ be a multiprojection of $U_k$ on a
$l$-dimensional subspace $W$, that we write as $mU_l=\widehat{mp}(W)U_k$.
It means that a matrix $M(mU_l)$ is obtained from the matrix $M(U_k)$ by
removal of columns that do not belong to $W$. We call a $k$-relation $U_k$
\emph{$mp$-symmetrical}, if $\widehat{mp}(W)U_k$ is homogenous (i.e each
line of $M(mU_l)$ has the same multiplicity) for any possible subspace
$W$. A $k$-partition $L_k$ is $mp$-symmetrical, if every its class is
$mp$-symmetrical.

We have described three necessary properties of an automorphic partition:
$s$-, $p$- and $mp$-symmetry, at that $p$-symmetry involves $s$-symmetry
(proposition \ref{p->s-symm}). A $k$-partition that possesses these three
symmetries we call \emph{regular} $k$-partition. A $k$-partition that is a
projection of a regular $k$-partition we call \emph{strongly-regular}. One
can see that regular and strongly-regular graphs satisfy corresponding
conditions. It is clear that strongly-regular partition is $pq$-stable.
Reverse statement is not correct, a counterexample is a $8$-point cubic
graph obtained from a cube in which two parallel edges $\{1,2\}$
and $\{3,4\}$, belonging to one cube face, are changed with edges
$\{1,3\}$ and $\{2,4\}$. This graph is point-transitive, its $2$-partition
$L_2$ (on edges and not edges) is assembling in $3$-partition $L_3$, but
$L_3$ is not $mp$-symmetrical. We will prove below the next

\begin{theorem}\label{pq->srp}
Let $k\geq 3$ and $L_k$ be a regular, $pq$-stable $k$-partition, then
$L_k$ is strongly-regular.

\end{theorem}


\section{Partition stabilization algorithm}

It is clear that any $k$-partition $L_k$ can be stabilized by
$pq$-stabilization to a $pq$-stable partition
$R_k=(\hat{p}\hat{q})^{\nu}L_k$, where $\nu$ is a number of iterations.
From theorem \ref{pq->srp} it follows that for $k\geq 3$ $R_k$ is
strongly-regular, if $L_k$ is regular. One can see that the algorithm of
the regularization of a $k$-partition follows immediately from its
definition and is polynomial. Concern of graph isomorphism it is of
interest whether exists a number $k$ for that $pq$-stabilization (of a
regular $k$-partition) leads to an automorphic $k$-partition or at least
to a strongly regular $k$-partition with non-trivial automorphism group.
If such number $k$ exists, then the graph isomorphism problem is
polynomial (because of complexity equivalence, considered above). We show
below that corresponding $k$ exists and is equal to $3$.


\section{Automorphic $k$-partitions}

Let $L_k$ be an automorphic $k$-partition, i.e. $L_k=Orb_k(Aut(L_k))$,
then we have the next

\begin{theorem}\label{L_(k+1)-auto}
Let $L_{(k+i)}=\hat{q}^iL_k$, $i\in [1,n-k]$, then $L_{(k+i)}$ is
automorphic.

\end{theorem}
\proof
Since $L_k$ is automorphic, $\hat{p}^i\hat{q}^iL_k=L_k$ for any $i\in
[1,n-k]$. So $L_{(k+i)}$ is $k$-full and $pq$-stable (or $L_{(k+i)}=
\hat{p}\hat{q}L_{(k+i)}=\hat{q}\hat{p}L_{(k+i)}$) for any $i$. It follows
that $\hat{p}^{n-k-i}\hat{q}^{n-k-i}L_{(k+i)}=
L_{(k+i)}$.~$\Box$\bigskip

A permutation group $G(V)$ is called $k$-closed, if $Aut(Orb_k(G))=G$.

\begin{corollary}\label{k-closed}
Let $G$ be a permutation group, then it is $k$-closed group iff its
$n$-orbit is $k$-full.

\end{corollary}

Let $G$ be a $1$-closed group, then it is a cartesian product of symmetric
groups acting on a partition of $V$.


\section{Algebraic combinatorics of strongly regular $k$-partitions}

The purpose of this section is a proof of theorem \ref{pq->srp} and a
proof of the polynomial complexity of the algorithm, detecting graph
non-trivial automorphism (and therefore a proof of polynomiality of the
graph isomorphism problem).

As we wrote above the contemporary algebraic combinatorics is a theory of
strongly regular $2$-partitions (or one can say strongly regular color
digraphs) that have historically many other names. With certain
restriction with an additional condition one obtains distance regular
graphs.  But this theory cannot tell many about possible symmetries on
$k$-partitions, so in order to obtain such information one has to consider
$k$-partitions for $k>2$. The main difficulty of such undertaking is that
by $k>2$ one cannot apply especially good developed matrix theory. So in
order to find an approach to investigation of $k$-partitions we put a
question: what is the most important in representation of $k$-partition?
And an answer could be: of course, it is its coloring function.  Now we
begin a search for an appropriate coloring function.

\subsection{Level invariant transformation}

Let $A=\{a_1,\ldots,a_d\}$ and $B=\{b_1,\ldots,b_d\}$ be sets of colors,
$L_k$ be a $k$-partition with $d$ classes and $f:L_k\rightarrow
A,g:L_k\rightarrow B$ be bijections. Let $T$ be a transformation of $A$ to
$B$: $B=TA$, so that $T$ is also a transformation of $f$ to $g$: $g=Tf$.
Such transformations maintain $L_k$ or level surfaces of $f$. We call $T$
a level invariant transformation for function $f$. Such transformations,
applied to a coloring function of $k$ variables $f(x_1,\ldots,x_k)$, give
different possibility for algebraic approach to investigation of
$k$-partitions. Here it will be of interest for us two level invariant
transformations. One of them is a polynomial $T(a)=P_{d-1}(a)$ of degree
$d-1$ that is defined by the next system of linear equations:

\begin{equation}\label{lev-inv-trans-p}
b_i=x_0+x_1a_i+\ldots+x_{d-1}a_i^{d-1}, i\in [1,d].
\end{equation}

And the second is a matrix $T$ that transforms the vector $\vec{a}=\langle
a_1,\ldots,a_d\rangle$ to the vector $\vec{b}=\langle
b_1,\ldots,b_d\rangle$:

\begin{equation}\label{lev-inv-trans-l}
\vec{b}=T\vec{a}
\end{equation}

We say that two functions $f$ and $g$ are equivalent $f\sim g$, if they
have the same level surfaces: $L_k(f)=L_k(g)$.

\subsection{Non-linear number coloring function}

Let $L_k$ be a $pq$-stable $k$-partition and
$\sigma_k\equiv\|\sigma_{i_1\ldots i_k}\|$ be an associated coloring tensor
on $L_k$. Let $L_{k+1}=\hat{q}L_k$ and
$\sigma_{k+1}\equiv\|\sigma_{i_1\ldots i_{k+1}}\|$ be a tensor
associated with $L_{k+1}$. We can represent the tensor $\sigma_{k+1}$
through the tensor $\sigma_k$ as:

\begin{equation}\label{sigma_(k+1)-nl}
\sigma_{i_1\ldots i_{k+1}}=\sigma_{i_1\ldots i_k}
\prod_{\langle j_1\ldots j_{k-1}\rangle
\in\hat{p}\langle i_1\ldots i_k\rangle}\sigma_{j_1\ldots
j_{k-1}i_{k+1}}.
\end{equation}

Since $L_k$ is a projection of $L_{k+1}$, we can represent the tensor on
$L_k$ through the tensor on $L_{k+1}$ as:

\begin{equation}\label{sigma_k-pq}
   \sigma_{i_1\ldots i_k}
   =\sum_{\alpha\in C_{k+1}}\sigma_{i_1\ldots i_k*}(\alpha),
\end{equation}
where $C_{k+1}=\{\sigma_{i_1\ldots i_ki_{k+1}}\}$ is a set of colors on
$L_{k+1}$ and $\langle i_1\ldots i_k*\rangle(\alpha)$ is a
representative $(k+1)$-tuple of a class $\alpha$. If there exists no such
representative of a class $\alpha$, then $\sigma_{i_1\ldots
i_k*}(\alpha)=0$. The factor $\sigma_{i_1\ldots i_k}$ in
(\ref{sigma_(k+1)-nl}) does not change the structure of the product and
can be omitted. Then using described above the level invariant
transformation (\ref{lev-inv-trans-p}) we find an equation on a tensor of
$pq$-stable $k$-partition in form:

\begin{equation}\label{k-eq-prod-pq}
\sum_{\alpha\in C_{k+1}}
\prod_{\langle j_1\ldots j_{k-1}\rangle
\in\hat{p}\langle i_1\ldots i_k\rangle}\sigma_{j_1\ldots
j_{k-1}*}(\alpha)=\emph{P}_{d-1}(\sigma_{i_1\ldots i_k}).
\end{equation}

For a strongly regular $k$-partition $L_k$ (that is $mp$-symmetrical) we
can rewrite equality (\ref{sigma_k-pq}) as:

\begin{equation}\label{sigma_k-srp}
\sigma_{i_1\ldots i_k}=
\sum_{\alpha\in C_{k+1}}\sigma_{i_1\ldots i_k,*}(\alpha)r_{i_1\ldots
i_k*}(\alpha)=\sum_{l=1,n}\sigma_{i_1\ldots i_kl},
\end{equation}
where $r_{i_1\ldots i_k*}(\alpha)$ is a multiplicity of a $k$-tuple
$\langle i_1\ldots i_k\rangle$ in a multiprojection of the class of color
$\alpha$ of $L_{k+1}$ by removing the latter column in a matrix of the
class $\alpha$.

Thus for a strongly regular $k$-partition $L_k$ the equation
(\ref{k-eq-prod-pq}) takes a form:

\begin{equation}\label{k-eq-prod}
\sum_l\prod_{\langle j_1\ldots j_{k-1}\rangle\in\hat{p}\langle i_1\ldots
i_k\rangle}\sigma_{j_1\ldots j_{k-1}l}=\emph{P}_{(d-1)}(\sigma_{i_1\ldots
i_k}),
\end{equation}
where we assume that elements with equal indices are zero and $l\neq
i_1,\ldots,i_k$ (\emph{this condition is implied in corresponding sums
below}).

\subsection{Linear number coloring function}

We can represent the tensor $\sigma_{k+1}$ through the tensor $\sigma_k$
also as:

\begin{equation}\label{k-eq-sum}
\sigma_{i_1\ldots i_{k+1}}=\sigma_{i_1\ldots i_k}x_0+
\sum_{\langle j_1\ldots j_{k-1}\rangle
\in\hat{p}\langle i_1\ldots i_k\rangle}\sigma_{j_1\ldots
j_{k-1}i_{k+1}}x_{\nu_{\langle j_1\ldots j_{k-1}\rangle}},
\end{equation}
where $\nu_{\langle j_1\ldots j_{k-1}\rangle}=m$, $m\in [1,k]$ is the
index of coordinate in $\{i_1,\ldots, i_k\}\setminus \{j_1,\ldots,
j_{k-1}\}$ and $x_0,x_1,\ldots,x_k$ are free parameters.

Using this representation and linear transformation
(\ref{lev-inv-trans-l}), we obtain an equation on coloring tensor of
strongly regular $k$-partition $L_k$ in a form:

\begin{equation}\label{k-eq-l}
\sum_l\sum_{\langle j_1\ldots j_{k-1}\rangle\in\hat{p}\langle i_1\ldots
i_k\rangle}\sigma_{j_1\ldots j_{k-1}l}x_{\nu_{\langle j_1\ldots
j_{k-1}\rangle}}=T^{j_1\ldots j_k}_{i_1\ldots i_k}\sigma_{j_1\ldots j_k},
\end{equation}
where the summand $\sigma_{i_1\ldots i_k}x_0$ is omitted and the right part
of equation is a conventional sum by $j$-indices that represents a
coloring of~$L_k$.

We can rewrite (\ref{k-eq-l}) in the equivalence form as:

\begin{equation}\label{k-eq-l-equiv}
\sum_l\sum_{\langle j_1\ldots j_{k-1}\rangle\in\hat{p}\langle i_1\ldots
i_k\rangle}\sigma_{j_1\ldots j_{k-1}l}x_{\nu_{\langle j_1\ldots
j_{k-1}\rangle}}\sim \sigma_{i_1\ldots i_k}.
\end{equation}

\subsection{Projective convolution of tensors}

Now we introduce an operation on tensors that we applied in left part of
equality (\ref{k-eq-prod}).

Let $\{\sigma_k^j:j=[1,k]\}$ be a set of $k$ tensors  of a rank $k$, then
we introduce a convolution operation:

\begin{equation}\label{ten-inv-trans}
\sum_l\prod_{j=1,k}\sigma^j_{\alpha(j)l}\equiv
\sigma_k^1\diamond\ldots\diamond\sigma_k^k,
\end{equation}
where $\alpha(j)=\langle i_1\ldots i_{j-1}i_{j+1}\ldots i_k\rangle$

For a tensor of the rank 2 it is the conventional matrix product.

\subsection{(0,1)-Tensor coloring function}

Consider the coloring of $k$-partition $L_k$, $|L_k|=d$, through $d$
$(0,1)$-tensors $a_k^\alpha$ with elements $a_{i_1\ldots i_k}^\alpha\in
\{0.1\}$, $\alpha\in [1,d]$, $i_1,\ldots,i_k\in [1,n]$. The value of
$a_{i_1\ldots i_k}^\alpha$ is 1 if a $k$-tuple $\langle i_1\ldots
i_k\rangle$ belongs to the class $U_k(\alpha)$ of $L_k$, else
$a_{i_1\ldots i_k}^\alpha=0$. So these $(0,1)$-tensors are linear
independent and any linear combination of them is a coloring of $L_k$.

Let $L_k$ be strongly regular, then we obtain an equation:

\begin{equation}\label{0.1-a_k->a_k)}
a_k^{i_1}\diamond\ldots\diamond a_k^{i_k}=
\sum_{\alpha=1,d}\lambda_{i_1\ldots i_k}^\alpha a_k^\alpha
\end{equation}

It is a generalization of equation on associative scheme \cite{bannai}.
We try here only to show the possibility of $k$-partition representation
and do not develop corresponding theory. Now we consider some examples.

\subsection{Strongly regular $2$-partitions}

A strongly regular simple graph $\Gamma(n,m,\lambda,\mu)$ is a strongly
regular $2$-partition $L_2$ that consists of $d=2$ symmetrical classes
($\langle v_1v_2\rangle L_2\langle v_2v_1\rangle$ for any $v_1,v_2\in V$).
If $L_3=\hat{q}L_2$, then parameters of $\Gamma$ are $n=|V|$, multiplicity
$m$ of a point in a class of $L_2$ (that represent edges in the graph) and
multiplicities $\lambda,\mu$ of pairs from two classes of $L_2$ in
corresponding two classes of $L_3$ (that represent triangles with 3 and 1
edges in the graph).

It is known that an adjacency matrix $A$ of the graph $\Gamma$ satisfies
to equation

\begin{equation}\label{srg-eq}
A^2= mE+\lambda A+\mu \bar{A},
\end{equation}
where $E$ is the unity matrix, $\bar{A}=I-A$ and $I$ is the $1$-matrix
($\{I_{ij}\}=\{1\}$).

The equation (\ref{srg-eq}) follows also from (\ref{k-eq-prod}). For
$k=2$ and $d=2$ we have a system of equations (for $i\neq j$)

\begin{equation}\label{2-2}
\sum_l \sigma_{il}\sigma_{jl}= x+y\sigma_{ij}
\end{equation}

Let $\sigma_{ij}\in \{\mu_0,\lambda_0\}$ and the tensor
$\sum_l\sigma_{il}\sigma_{jl}$ has values $\lambda$ and $\mu$ that color
classes $\lambda_0$ and $\mu_0$ of $L_2$ correspondingly, then we obtain
the system of equations of coloring transformation:

\begin{eqnarray*}
\lambda= x+y\cdot \lambda_0\\
\mu=x+y\cdot \mu_0.
\end{eqnarray*}

So the system (\ref{2-2}) takes a form:

$$
\sum_l\sigma_{il}\sigma_{jl}=
\frac{\mu\lambda_0-\mu_0\lambda}{\lambda_0-\mu_0}+
\frac{\lambda-\mu}{\lambda_0-\mu_0}\sigma_{ij}.
$$

For the 0,1-tensor $\sigma_{ij}$ ($\mu_0=0$ and $\lambda_0=1$) we obtain

\begin{equation}\label{srg-eq-ij}
\sum_l\sigma_{il}\sigma_{jl}=\mu+(\lambda-\mu)\sigma_{ij}.
\end{equation}

Or using, matrix notation ($\|\sigma_{ij}\|=A$)

$$
A^2-mE=\mu I+(\lambda-\mu)A=\lambda A+\mu \bar{A}.
$$

\subsection{$pq$-Stable regular $3$-partition}

Let $L_3$ be a $pq$-stable regular $3$-partition and $\sigma_{ijk}$ be its
coloring tensor. Let $L_4=\hat{q}L_3$, then there exists a coloring
$\sigma_{ijkl}$ of $L_4$ that can be represented accordingly to
(\ref{sigma_(k+1)-nl}) as:

\begin{equation}\label{3->4}
\sigma_{ijkl}=\sigma_{ijk}\sigma_{ijl}\sigma_{ikl}\sigma_{jkl}\,.
\end{equation}

Let $L_2=\hat{p}L_3$ and $\sigma_{ij}$ be a tensor on $L_2$, representing
$\sigma_{ijk}$ as:

\begin{equation}\label{2->3-nl}
\sigma_{ijk}=\sigma_{ij}\sigma_{ik}\sigma_{jk}\,,
\end{equation}
then

\begin{equation}\label{2->4}
\sigma_{ijkl}=\sigma_{ij}^2\sigma_{ik}^2\sigma_{jk}^2
\sigma_{il}^2\sigma_{jl}^2\sigma_{kl}^2.
\end{equation}

We consider here the $mp$-symmetry of $L_4$. In this case a power and
first three factors of equality (\ref{2->4}) are not of principal, so we
can rewrite that as:

\begin{equation}\label{2->4,1}
\sigma_{ijkl}=\sigma_{il}\sigma_{jl}\sigma_{kl}.
\end{equation}

Let $\langle ijk\rangle,\langle i'j'k'\rangle\in U_3\in L_3$ and $\langle
ijkl\rangle\in U_4\in L_4$. We will show that the multiplicity of
$3$-tuples $\langle ijk\rangle$ and $\langle i'j'k'\rangle$ in
$4$-relation $U_4$ are equal. Consider sums $S(\langle ijk\rangle)=
\sum_l\sigma_{il}\sigma_{jl}\sigma_{kl}$ and $S(\langle i'j'k'\rangle)=
\sum_l\sigma_{i'l}\sigma_{j'l}\sigma_{k'l}$. Since $\sum_m\sigma_{ml}=
\sigma_l$ is a coloring of $L_1=\hat{p}L_2$, then we find that
$\sum_kS(\langle ijk\rangle)=\sum_{k'}S(\langle i'j'k'\rangle)$,
$\sum_jS(\langle ijk\rangle)=\sum_{j'}S(\langle i'j'k'\rangle)$ and
$\sum_iS(\langle ijk\rangle)=\sum_{i'}S(\langle i'j'k'\rangle)$. These
equalities are valid by different coloring tensor of $L_2$, hence
$S(\langle ijk\rangle)= S(\langle i'j'k'\rangle)$ also for different
coloring of $L_2$. It proves theorem \ref{pq->srp}.

\subsection{Strongly regular $3$-partitions}

Here we show that

\begin{theorem}\label{srg-3-p}
A non-trivial strongly regular $3$-partition contains non-trivial
automorphism.

\end{theorem}
and therefore prove the polynomiality of graph isomorphism problem.

Let $L_3$ be a strongly regular $3$-partition and $\sigma_{ijk}$ be its
coloring tensor, then the equivalence (\ref{k-eq-l-equiv}) takes a form:

\begin{equation}\label{3->3+}
   \sum_l(\sigma_{ijl}z+\sigma_{ikl}y+\sigma_{jkl}x)\sim\sigma_{ijk}.
\end{equation}

Let $L_3$ be $2$-full, then we can represent the coloring tensor of $L_3$
through some coloring tensor of $L_2=\hat{p}L_3$, so that

\begin{equation}\label{2->3-l}
\sigma_{ijk}=\sigma_{ik}p+\sigma_{jk}q,
\end{equation}
where $p,q$ are also free parameters.

By substitute the right part of (\ref{2->3-l}) for the tensor of $L_3$ in
(\ref{3->3+}) we obtain an equivalence:

\begin{equation}\label{S3}
\sum_l(\sigma_{il}x+\sigma_{jl}y+\sigma_{kl}z)\sim
\sigma_{ik}p+\sigma_{jk}q=\sigma_{ijk}
\end{equation}
(here x,y,z are new parameters).

Now we consider what for equations follow from equivalence (\ref{S3}).
Let tensor $\|\sigma_{ij}\|$ (and correspondingly tensor
$\|\sigma_{ijk}\|$) has non-trivial automorphism $\phi$, then
$$\sigma_{\phi(i)\phi(j)\phi(k)}=\sigma_{ijk}$$ and
$$\sigma_{\phi(i)\phi(l)}x+\sigma_{\phi(j)\phi(l)}y+
\sigma_{\phi(k)\phi(l)}z=\sigma_{il}x+\sigma_{jl}y+\sigma_{kl}z.$$ So we
have a bijection between summands of sums that represents $\sigma_{ijk}$
and $\sigma_{\phi(i)\phi(j)\phi(k)}$ in equivalent coloring. It gives a
possibility to reduce in (\ref{S3}) the number of independent variables
$\sigma_{ij}$ and at the same time the number of independent related by
equivalence lines (for different $3$-tuples $\langle ijk\rangle$), by
substitution anywhere in (\ref{S3}) $\sigma_{ij}$ for
$\sigma_{\phi(i)\phi(j)}$ ($\phi\in Aut(\|\sigma_{ij}\|)$, $i,j\in [1,n]$).

Let now $2$-partition $L_2=L(\|\sigma_{ij}\|)$ and $3$-partition
$L_3=\hat{q}L_2=L(\|\sigma_{ijk}\|)$ be faithful strongly regular, then
there exist $3$-tuples $\langle ijk\rangle$ and $\langle i'j'k'\rangle$
that belong to the same class of $L_3$ and are connected with no
automorphism. Then $\sigma_{ijk}=\sigma_{i'j'k'}$ and

\begin{equation}\label{f-eq-1}
\sum_l(\sigma_{il}x+\sigma_{jl}y+\sigma_{kl}z)=
\sum_l(\sigma_{il}x+\sigma_{jl}y+\sigma_{kl}z),
\end{equation}
Since $x,y,z$ are free parameters then from (\ref{f-eq-1}) it follows
equalities of three subsums:
\begin{equation}\label{f-eq-1x}
\sum_{l\neq i,j,k}\sigma_{il}=\sum_{l\neq i',j',k'}\sigma_{i'l},
\end{equation}
\begin{equation}\label{f-eq-1y}
\sum_{l\neq i,j,k}\sigma_{jl}=\sum_{l\neq i',j',k'}\sigma_{j'l},
\end{equation}
\begin{equation}\label{f-eq-1z}
\sum_{l\neq i,j,k}\sigma_{kl}=\sum_{l\neq i',j',k'}\sigma_{k'l}.
\end{equation}
Because of s-symmetry this three systems of subsums equalities (for
different pairs of $3$-tuples) are equal. Thus it is sufficient to consider
the system $S_x$ given by expression (\ref{f-eq-1x}) and choose only such
equations in this system that are independent by automorphisms and by
transitivity. This system of equations has solution if the number of
equations $|S_x|$ is less than number of variables, because the equalities
in $S_x$ are independent. From here immediately follows that, in the case
of strongly regular $3$-partition $L_3$, the system $S_x$ can be solved if
$Aut(L_3)$ is enough rich on automorphisms, because, when $Aut(L_3)$ is
trivial, the number of equations in $S_x$ is $O(n^3)$ and the numner of
variables is $O(n^2)$.  This proves theorem \ref{srg-3-p}.


\section*{Conclusion}

In given solution of graph isomorphism problem were used symmetry
properties of $k$-orbits. Other texts of author connected with
consideration of $k$-orbits one can find in "www.arxiv.org". Those texts
are not mistake free, but they contains new original ideas and a
direction of investigation, and therefore could be of interest.
Author hopes that investigation of symmetry properties of $k$-orbits can
bring new ideas for simplifying of simple finite group classification.


\section*{Acknowledgements}

I would like to express my thanks to Dr. M.Tabachnikov (Kharkov, Ukraine)
for the proposal (in 1983) to solve the graph isomorphism problem, to Dr.
V.Grinberg (Kharkov-USA) for productive contacts, to Dr. M.Klin
(Moscow-Israel) for acquainting with contemporary achievements in this
part of combinatorics and to others, who helped me on the way.


\end{document}